\documentclass[10pt,a4paper,leqno]{article}
\usepackage[french]{babel}
\usepackage{amsmath,amstext, amsthm, amssymb}
\usepackage{amsmath}
\usepackage{mathrsfs}
\usepackage{amsmath,amstext, amsthm, amssymb}

\def\og{\leavevmode\raise.3ex\hbox{$\scriptscriptstyle\langle\!\langle$~}}
\def\fg{\leavevmode\raise.3ex\hbox{~$\!\scriptscriptstyle\,\rangle\!\rangle$}}

\font \tenmat = msbm10 \font \sevenmat = msbm7 \font \fivemat =
msbm5 \newfam \matfam \textfont \matfam = \tenmat \scriptfont
\matfam = \sevenmat \scriptscriptfont \matfam = \fivemat

\def\RR{{\mathbb R}} 
\def\CC{{\mathbb C}}

\def\cH{{\mathcal H}} 
\def\cM{{\mathcal M}} 
\def\li{{\rm Li}_2} 

\newtheorem{theorem}{Theorem}[section]

\newtheorem{e-proposition}[theorem]{Proposition}

\newtheorem{e-definition}[theorem]{Definition\rm}

\title{On the Modular Behaviour of the Infinite Product $(1-x)(1-xq)(1-xq^2)(1-xq^3)\cdots$.}
\author{Changgui ZHANG}
\date{}

\begin{document}

\maketitle
 
\begin{abstract}
 Let $q=e^{2\pi i\tau}$, $\Im\tau>0$, $x=e^{2\pi i\xi}\in\CC$ and $(x;q)_\infty=\prod_{n\ge 0}(1-xq^n)$. Let $(q,x)\mapsto(q^*,\iota_q x)$ be the classical modular substitution given by $q^*=e^{-2\pi i/\tau}$ and  $\iota_q x=e^{2\pi i\xi/{\tau}}$. The main goal of this Note is to study the ``modular behaviour'' of the infinite product $(x;q)_\infty$, this means, to compare the function defined by $(x;q)_\infty$ with that given by $(\iota_q x;q^*)_\infty$. Inspired by the work  \cite{St} of Stieltjes on some semi-convergent series, we are led to a ``closed'' analytic formula for $(x;q)_\infty$ by means of the dilogarithm combined with a Laplace type integral that admits a divergent series  as Taylor expansion at $\log q=0$. Thus, we can obtain an expression linking $(x;q)_\infty$ to its modular transform $(\iota_qx;q^*)_\infty$ and which contains, in essence, the modular formulae known for Dedekind's eta function, Jacobi theta function and also for certain Lambert series. Among other applications, one can  remark  that our results allow to obtain a Ramanujan's asymptotic formula about $(x;q)_\infty$ for $q\to 1$  \cite[Entry 6, p. 265 \& Entry 6', p. 268]{Be1}.

\bigskip
\bigskip

{\bf Sur le comportement modulaire du produit infini $(1-x)(1-xq)(1-xq^2)(1-xq^3)\cdots$.}
Soit $q=e^{2\pi i\tau}$, $\Im\tau>0$, $x=e^{2\pi i\xi}\in\CC$ et $(x;q)_\infty=\prod_{n\ge 0}(1-xq^n)$. Soit $(q,x)\mapsto (q^*,\iota_q x)$ la classique substitution modulaire donn\'ee par $q^*=e^{-2\pi i/\tau}$ et  $\iota_q x=e^{2\pi i\xi/{\tau}}$. Le principal but de la pr\'esente Note est d'\'etudier le \og comportemant modulaire \fg\ du produit infini $(x;q)_\infty$, c'est-\`a-dire, de comparer la fonction d\'efinie par $(x;q)_\infty$ \`a celle par $(\iota_q x;q^*)_\infty$. Inspir\'e du travail \cite{St} de Stieltjes sur des s\'eries semi-convergentes, nous somme parvenus \`a une formule analytique \og explicite\fg\ pour $\log(x;q)_\infty$ au moyen du dilogarithme complet\'e par une int\'egrale du type Laplace, cette derni\`ere ayant une s\'erie divergente dans son d\'eveloppement taylorien en $\log q=0$. Ceci nous permet d'obtenir une expression reliant $(x;q)_\infty$ \`a sa transform\'ee modulaire $(\iota_qx;q^*)_\infty$ et qui contient essentiellement les formules modulaires connues pour la fonction eta de Dedekind, la fonction theta de Jacobi et aussi pour certaines s\'eries de Lambert. Parmi d'autres applications, on remarquera que nos r\'esultats permettent  \`a obtenir  une formule asymptotiquede de Ramanujan sur  $(x;q)_\infty$  pour $q\to 1$ \cite[Entry 6, p. 265 \& Entry 6', p. 268]{Be1}.

\end{abstract}

\section*{Version fran\c{c}aise abr\'eg\'ee}

Soit $q=e^{2\pi i\tau}$, $\Im\tau>0$, $x=e^{2\pi i\xi}\in\CC$ et consid\'erons le produit infini figurant sur le titre de la Note, not\'e d\'esormais $(x;q)_\infty$. L'\'etude de ce dernier peut remonter au temps d'Euler \cite[Chap. XVI]{Eu}, \`a qui l'on doit l'identit\'e remarquable \eqref{equation:Euler} rappel\'ee plus loin.
Celle-ci justifie que $((q-1)x;q)_\infty$ est un $q$-analogue de la fonction exponentielle, avec $(q;q)_n/(1-q)^n\to n!$ pour $q\to 1$. Le produit infini $(x;q)_\infty$ occupe ainsi une place centrale dans le monde des \og $q$-s\'eries \fg, y compris la th\'eorie des fonctions ellptiques \cite{WW} ou celle des \'equations aux $q$-diff\'erences  \cite{Zh}, \cite{DZ}, \cite{RSZ}.

Soit la transformation  \og modulaire\fg\ $\cM:\ (q,x)\mapsto(q^*,\iota_qx)=\bigl(e^{-2\pi i/\tau},e^{2\pi i\xi/{\tau}}\bigr)$; afin de comparer la fonction $(x;q)_\infty$ \`a sa transform\'ee modulaire $(\iota_qx;q^*)_\infty$, commen\c cons par les observations suivantes.

\begin{description}
 \item (i) Dans le dernier paragraphe \cite[p. 252-258]{St} de l'article de Stieltjes intitul\'e \og Recherches sur quelques s\'eries semi-convergentes\fg,  une int\'egrale \og singuli\`ere \`a la Cauchy\fg, {\it i.e}, avec des valeurs principales, a \'et\'e donn\'ee pour repr\'esenter une s\'erie de Lambert $P(a)$ et a permis de lier $P(a)$ \`a $P(1/a)$.
\item (ii) Le produit infini $(q;q)_\infty$,  li\'e \`a la fonction eta de Dedekind ou encore \`a la fonction $\tau$ de Ramanujan, est muni d'une formule modulaire  \cite[(44), p. 154]{Se} et celui-ci est aussi la valeur de $(x;q)_\infty$ au point $x=q$.
\item (iii) La fonction theta de Jacobi v\'erifie une relation modulaire et elle peut s'\'ecrire comme le produit de $(q;q)_\infty$ par le facteur $(\sqrt qx;q)_\infty(\sqrt q/x;q)_\infty$, lequel est invariant par $x\mapsto1/x$.
\end{description}

Cela \'etant, nous sommes amen\'es \`a la question de savoir si la m\'ethode (i) de Stieltjes sera susceptible d'\^etre g\'en\'eralis\'ee \`a la fonction $(x;q)_\infty$ ou $\log(x;q)_\infty$ pour obtenir une relation, de la forme 
$$(x;q)_\infty=K(q,x)(\iota_qx;q^*)_\infty,\eqno(K)
$$ 
v\'erifiant les propri\'et\'es ci-dessous. 

\begin{description}
 \item (iv) Le facteur multiplicatif $K(x,q)$ sera donn\'e \og autant explicite\fg\ que possible.
\item (v) Losque $x=q$ ou que les facteurs $(\sqrt qx;q)_\infty$ et $(\sqrt q/x;q)_\infty$ sont mis ensemble, la relation (K) voulue ci-dessus nous conduira aux formules modulaires rappel\'ees  dans (ii) et (iii) ci-dessus respectivement.
\end{description}

Le principal but de la Note est de montrer que la m\'ethode de Stieltjes \'evoqu\'ee ci-dessus r\'epond \`a notre probleme et que le facteur $K(q,x)$ demand\'e dans la relation (K) sera donn\'e au moyen du dilogarithme avec une int\'egrale de Laplace, ce qui prouve \'egalement l'aspect \og semi-convergent\fg\ de $(x;q)_\infty$ pour $q\to 1$, $\vert q\vert<1$; voir Theorem \ref{theorem:xqmain} et la relation correspondante \eqref{equation:xqmain} ou les formulations \eqref{equation:K(q,x)}-\eqref{equation:K-}.
Par ailleurs, vu les relations \eqref{equation:symetricG*P}-\eqref{equation:Landen}, nous en d\'eduisons ais\'ement la formule modulaire connue de la fonction $\eta$ de Dedekind et aussi celle de la fonction $\theta$ de Jacobi. En consid\'erant la d\'eriv\'ee logarithmique par rapport \`a la variable $x$ pour chacun des c\^ot\'es de \eqref{equation:xqmain}, nous en obtenons des formulas qui relient chaque s\'erie de Lambert (g\'en\'eralis\'ee) \`a la s\'erie obtenue apr\`es la substitution $(q,x)\mapsto\cM(q,x)$.

\section{Introduction}

There are more and more studies on $q$-series and related topics, not only in traditional themes, but also in more recent branches, such as quantum physics, random matrices. A first non-trivial example of $q$-series may be the infinite product $(q;q)_\infty:=\prod_{n\ge 1}(1-q^n)$, that is considered in Euler \cite[Chap. XVI] {Eu} and then is revisited by many of his successors, particularly intensively by Hardy and Ramanujan \cite[p. 238-241; p. 276-309; p. 310-321]{Ha}. Beautiful formulae are numerous and motivations are often various: elliptic and modular functions theory, number and partition theory, orthogonal polynomials theory, {\it etc} $\dots$. Concerning this wonderful history, one may think of Euler's pentagonal number theorem \cite[p. 30] {An}, Jacobi's triple product identity \cite[\S10.4, p. 496-501]{AAR}, Dedekind  modular eta function \cite[(44), p. 154]{Se}, \cite[Chapter 3] {Ap}, to quote only some examples of important masterpieces.

However, the  infinite product  $(x;q)_\infty:=\prod_{n\ge 0}(1-xq^n)$, also already appearing as initial model in the same work \cite[Chap. XVI]{Eu} of Euler, receives less attention although it always plays a remarkable role in all above-mentioned subjects:  thanks to Euler, one knows the following remarkable formula:
\begin{equation}\label{equation:Euler}
(x;q)_\infty=\sum_{n\ge 0}\frac{q^{n(n-1)/2}}{(q;q)_n}\,(-x)^n\qquad \bigl((q;q)_n=(q;q)_\infty/(q^{n+1};q)_\infty\bigr)\,,
\end{equation}
which says in what manner $((q-1)x;q)_\infty$ may be viewed as a $q$-analog of the standard exponential function; see \cite{Zh}, \cite{DZ} and \cite{RSZ} for a point of view from the analytic theory of $q$-difference equations. Moreover, $(x;q)_\infty$ is immediately linked with several constants in the elliptic integral theory, Gauss' binomial formula \cite[\S 10.2, p. 487-491]{AAR}, generalized Lambert series,  {\it etc} $\dots$. Indeed, rapidly the situation may become more complicated and, what is really important, the modular relation does not occur for generic values of $x$.

In this Note, we shall point out how, up to an explicit part, the function defined by the product $(x;q)_\infty:=(1-x)(1-xq)(1-xq^2)(1-xq^3)...$ can be seen somewhat modular. This {\it non-modular part} can be considered as being represented by a divergent but Borel-summable power series on variable $\log q$ near zero, that is, when $q$ tends toward the unit value. These results, subjects of Theorems \ref{theorem:main} and \ref{theorem:xqmain}, give rise to new and unified approaches to treat Jacobi theta function, Lambert series, $\dots$.

\bigskip
The Note is organized as follows. In Section \ref{sec:main}, we will give some integral representations about the real-valued function $\log(x;q)_\infty$ for $q\in(0,1)$ and $x\in(0,1)$. In Section \ref{section:main}, we will extend the above results to the complex case and therefore find the modular type formula for $(x;q)_\infty$. The main formulas are \eqref{equation:main} and \eqref{equation:xqmain}.

\section{Integral representations of $\log(x;q)_\infty$}\label{sec:main}
In this section, we will suppose that $q\in(0,1)$ and $x\in(0,1)$, so that the infinite product $(x;q)_\infty$ takes values in $(0,1)$. By applying an idea of Stieltjes \cite[p. 252-258]{St} to the real-valued function $\log(x;q)_\infty$, we firstly find the following statement.

\begin{theorem}\label{theorem:main}
Let $q=e^{-2\pi\alpha}$, $x=e^{-2\pi(1+\nu)\alpha}$ and suppose $\alpha>0$ and $\nu>-1$. The following relation holds:
\begin{align}\label{equation:main}
\log(x;q)_\infty=&-\frac{\pi}{12\alpha}+\log\frac{\sqrt{2\pi}}{\Gamma(\nu+1)}+\frac \pi{12}\,\alpha-\bigl(\nu+\frac12\bigr)\log\frac{1-e^{-2\pi\nu \alpha}}\nu
+\int_0^{\nu }\bigl(\frac {2\pi\alpha t}{e^{2\pi\alpha t}-1}-1\bigr)\,dt+M(\alpha,\nu),
\end{align}
where
\begin{equation}\label{equation:M}
M(\alpha,\nu)=-\sum_{n=1}^\infty\frac{\cos2n\pi\nu}{n(e^{2n\pi/\alpha}-1)}-\frac2\pi\,{\mathcal PV}\int_0^\infty\sum_{n=1}^\infty\frac{\sin2n\nu \pi t}{n(e^{2n\pi t/\alpha}-1)}\,\frac{dt}{1-t^2} \,.
\end{equation}
\end{theorem}

In the above, $\displaystyle {\mathcal PV}\int$ stands for the principal value of a singular integral in the Cauchy's sense; see \cite[\S 6.23, p. 117]{WW}. 
If we let
\begin{equation}\label{equation:B}
B(t)=\frac1{e^{2\pi t}-1}-\frac1{2\pi t}+\frac12\,,
\end{equation}
it is easy to see Theorem \ref{theorem:main} can be stated as follows.

\begin{theorem}\label{theorem:mainbis}
Let $q$, $x$, $\alpha$, $\nu$ and $M(\alpha,\nu)$ be as given in Theorem \ref{theorem:main}. Then the following relation holds:
\begin{align*}\label{equation:mainbis}
\log(x;q)_\infty=&-\frac{\pi}{12\alpha}-\bigl(\nu+\frac12\bigr)\,\log2\pi\alpha+\log\frac{\sqrt{2\pi}}{\Gamma(\nu+1)}+\frac\pi2\,(\nu+1)\,\nu\alpha\cr
&+\frac \pi{12}\,\alpha+2\pi\alpha\int_0^\nu\bigl(t-\nu-\frac12\bigr)\,B(\alpha t)\,dt +M(\alpha,\nu),
\end{align*}
where $B$ is defined in \eqref{equation:B}.
\end{theorem}

As usual, let $\li$ denote the dilogarithm function; recall $\li$ can be defined as follows \cite[(2.6.1-2), p. 102]{AAR}:
\begin{equation}\label{equation:dilog}
\li (x)=-\int_0^x\log(1- t)\,\frac{dt}{t}=\sum_{n= 0}^{\infty}\frac{x^{n+1}}{(n+1)^2}\,.
\end{equation}

\begin{theorem}\label{theorem:mainter}
The following relation holds for any $q\in(0,1)$ and $x\in(0,1)$:
\begin{align}\label{equation:mainter}
\log(x;q)_\infty=\frac1{\log q}\,{\li(x)}+\log\sqrt{1-x}-\frac{\log q}{24}-\int_0^\infty B(-\frac{\log q}{2\pi}\,t)\,x^{t}\,\frac{dt}t+M(-\frac{\log q}{2\pi}\,,\log_qx)\,,
\end{align}
where $B$ denotes the function given by \eqref{equation:B}.
\end{theorem}

We shall write the singular integral part in \eqref{equation:M} by means of contour integration in the complex plane.
Fix a real $r\in(0,1)$ and, for any positive integer $n$, let $C_{n,r}^-$ ($C_{n,r}^+$, resp.) be the half circle passing from $n-r$ to $n+r$ by the right (left, resp.) hand side. Let $\ell^{\mp}_r=(0,1-r)\cup\Bigl(\cup_{n\ge 1}\bigl(C_{n,r}^\mp\cup(n+r,n+1-r)\bigr)\Bigr)$ and define $P^\mp(\alpha, \nu)$ as follows:
\begin{equation}\label{equation:Pmp}
P^\mp(\alpha,\nu)=\int_{\ell_r^{\mp}}\frac{\sin\nu t}{e^{ t/\alpha}-1}\,\big(\cot\frac{t}2-\frac2{t}\bigr)\,\frac{dt}{t}\,,
\end{equation} 
where $\alpha>0$ and where $\nu$ may be an arbitrary  real number.

\begin{theorem}\label{theorem:MP}
Let $M$ be as in Theorem \ref{theorem:main} and let $P^-$ be as in \eqref{equation:Pmp}. For any $\nu\in\RR$ and $\alpha>0$, the following relation holds:
\begin{equation*}\label{equation:MP}
M(\alpha,\nu)=\log(e^{2 \pi\nu i-2\pi /\alpha};e^{-2\pi /\alpha})_\infty+P^-(\alpha,\nu)\,.
\end{equation*}
\end{theorem}

Consequently, the term $M$ appearing in Theorem \ref{theorem:main} can be considered as being an \emph{almost modular term} of $\log(x;q)_\infty$;  the correction term $P^-$ given by \eqref{equation:Pmp} will be called {\it disruptive factor or perturbation term}.

\section{Modular type expansion of $(x;q)_\infty$}\label{section:main}
From now on, we will work with complex variables $q$ and $x$, with $q=e^{2\pi i\tau}$, $\tau\in\cH$ and $x=e^{2\pi i\xi}\in\CC^*$.
As usual, $\log$ will stand for the principal branch of the logarithmic function over its Riemann surface  $\tilde\CC^*$. For any pair of real numbers $a<b$, we will define
\begin{equation}\label{equation:Sab}
S(a,b):=\{z\in\tilde\CC^*:\arg z\in(a,b)\}\,;
\end{equation}
therefore, the Poincar\'e's half-plane $\cH$ will be identified to $S(0,\pi)$ while the  broken plane $\CC\setminus(-\infty,0]$ will be seen as the subset $S(-\pi,\pi)\subset\tilde\CC^*$.

Firstly, let us introduce the following modified complex version of $P^-$ given by \eqref{equation:Pmp}: for any $d\in(-\pi,0)$, let
\begin{equation}\label{equation:Pd}
P^d(\tau,\xi)=\int_{0}^{\infty e^{id}}\frac{\sin\frac\xi\tau\, t}{e^{i t/\tau}-1}\,\big(\cot\frac{t}2-\frac2{t}\bigr)\,\frac{dt}{t}\,,
\end{equation} 
the path of integration being the half line starting from origin to infinity with argument $d$. Note that  $P^d$ is analytic over the domain $\Omega^d$ if we set
\begin{align*}\label{equation:Omegad}
\Omega^d=\cup_{\sigma\in(0,\pi)}\bigl(0,\infty e^{i(d+\sigma)}\bigr)\times\{\xi\in\CC:\bigl\vert\Im(\xi\,e^{-i\sigma})\bigr\vert<\sin\sigma\}\,.
\end{align*}

The family of functions $\{P^d\}_{d\in(-\pi,0)}$ given by \eqref{equation:Pd} gives rise to an analytical function over the domain \begin{equation}\label{equation:Omega-}
\Omega_-:=S(-\pi\,,\pi)\times\Bigl( \CC\setminus\bigl((-\infty,-1]\cup[1,\infty)\bigr)\Bigr)\subset\CC^2\,.
\end{equation}
Moreover, if we denote this function by $P_-(\tau,\xi)$, then the following relation holds for all $\alpha>0$ and $\xi\in\RR$:
\begin{equation}\label{equation:PP-}
 P_-(\alpha i,\xi\alpha i)=P^{-}(\alpha,\xi)\,.
\end{equation} On the other hand, if we take the arguments $d\in(0,\pi)$ instead of $d\in(-\pi,0)$ in \eqref{equation:Pd}, we can get an analytical function, say $P_+$, defined over 
$$
\Omega_+:=S(0\,,2\pi)\times\Bigl( \CC\setminus\bigl((-\infty,-1]\cup[1,\infty)\bigr)\Bigr)\,.
$$
Therefore, one gets the following Stokes relation related to the above Lapalce type integral \eqref{equation:Pd}.

\begin{theorem}\label{theorem:StokesP}
For any $\tau\in\cH$, the relation
\begin{equation*}\label{equation:StokesP+-} 
P_-(\tau ,\xi)-P_+(\tau ,\xi)={2i}\,\sum_{n=1}^\infty\frac{\sin\frac{2n\xi\pi}\tau}{n(e^{2n\pi i/\tau}-1)}
\end{equation*}
holds provided that $|\Im(\xi/\tau)|<-\Im(1/\tau)$.
\end{theorem}

Secondly, the integral term involving the function $B$ in formula \eqref{equation:mainter} of Theorem \ref{theorem:mainter} is related to the remainder term of the Stirling's formula on $\Gamma$-function. An elementary analysis leads to introducing the following function
\begin{equation}\label{equation:GGamma}
G(\tau,\xi)=-\log\Gamma(\frac\xi\tau+1)+\bigl(\frac\xi\tau+\frac12\bigr)\,\log\frac\xi\tau-\frac\xi\tau+\log\sqrt{2\pi}\,,
\end{equation}
for all pair $(\tau,\xi)\in U^+:=\bigl\{(\tau,\xi)\in\CC^*\times\CC^*: \xi/\tau\notin(-\infty,0]\bigr\}$. One can check the following relation:
\begin{equation}\label{equation:G+}
G(\tau,\xi)+G(\tau,-\xi)=\log(1-e^{\mp 2\pi i \xi/\tau})
\end{equation}
according to $\displaystyle \frac\xi\tau\in S(-\pi,0)$ or $S(0,\pi)$, respectively.

Therefore, we come back to the extension of Theorem \ref{theorem:mainter} into the complex plane. for both complex. 
We recall the notations $q^*=e^{-2\pi i/\tau}$ and $\iota_qx=e^{2\pi i\xi/\tau}$; we will write $x^*$ instead of $\iota_qx$.

\begin{theorem}\label{theorem:xqmain}
The following relation holds for any $\tau\in\cH$ and $\xi\in\CC\setminus\bigl((-\infty,-1]\cup[1,\infty)\bigr)$ such that $\xi/\tau\notin(-\infty,0]$:
\begin{equation}\label{equation:xqmain}
(x;q)_\infty=q^{-1/24}\,\sqrt{1-x}\,\,(x^*q^*;q^*)_\infty\,\exp\Bigl(\frac{\li(x)}{\log q}+G(\tau,\xi)+P(\tau,\xi)\Bigr)\,,
\end{equation}
where $\sqrt{1-x}$ stands for the principal branch of $e^{\frac12\log(1-x)}$, $\li$ denotes the dilogarithm recalled in \eqref{equation:dilog}, $G$ is given by \eqref{equation:GGamma} and where $P$ denotes the function $P_-$ defined over $\Omega_-$ as in \eqref{equation:Omega-}-\eqref{equation:PP-}.
\end{theorem}

If we denote by $G^*$ the anti-symetrization of $G$ given by
\begin{equation*}
G^*(\tau,\xi)=\frac12\,\bigl(G(\tau,\xi)-G(\tau,-\xi)\bigr)\,,
\end{equation*}
then, according to relation \eqref{equation:G+}, we may rewrite \eqref{equation:xqmain} as follows:
\begin{equation}\label{equation:K(q,x)}
(x;q)_\infty=K(q,x)(x^*;q^*)_\infty\,,
\end{equation}
where the factor $K(q,x)$ is given in the following manner:
\begin{equation}\label{equation:K}
K(q,x)=q^{-1/24}\,\sqrt{\frac{1-x}{1-x^*}}\,\exp\Bigl(\frac{\li(x)}{\log q}+G^*(\tau,\xi)+P(\tau,\xi)\Bigr)
\end{equation}
if $\xi\in \tau\cH$, and
\begin{equation}\label{equation:K-}
K(q,x)=q^{-1/24}\,\frac{\sqrt{(1-x)(1-1/x^*)}}{1-x^*}\,\exp\Bigl(\frac{\li(x)}{\log q}+G^*(\tau,\xi)+P(\tau,\xi)\Bigr)
\end{equation}
if $\xi\in-\tau\cH$, that is, if $\displaystyle\frac\xi\tau\in S(-\pi,0)$.

In the above,  $G^*$ and $P$ are odd functions on the variable $\xi$:
\begin{equation}\label{equation:symetricG*P}
G^*(\tau,-\xi)=-G^*(\tau,\xi),\quad
P(\tau,-\xi)=-P(\tau,\xi)\,;
\end{equation} 
$\li$ satisfies the so-called \emph{Landen's transformation} \cite[Theorem 2.6.1, p. 103]{AAR}:
\begin{equation}\label{equation:Landen}
\li(1-x)+\li(1-\frac1x)=-\frac12\,\bigl(\log x\bigr)^2\,.
\end{equation}

Finally, if we write $\vec\omega=(\omega_1,\omega_2)=(1,\tau)$ and denote by $\Gamma_2(z,\vec\omega)$ the Barnes' double Gamma function associated to the double period $\vec\omega$ (\cite{Ba}), then Thoerme \ref{theorem:xqmain} and Proposition 5 of \cite{Sh} imply that
\begin{align*}\label{equation:doubleGamma}
\frac{\Gamma_2(1+\tau-\xi,\vec\omega)}{\Gamma_2(\xi,\vec\omega)}=&\sqrt{i}\,\sqrt{1-x}\,\exp\Bigl(\frac{\pi i}{12\tau}+\frac{\pi i}{2}\bigl(\frac{\xi^2}{\tau}-(1+\frac1\tau)\xi\bigr)+\frac{\li(x)}{\log q}+G(\tau,\xi)+P(\tau,\xi)\Bigr)\cr
=&\,\sqrt{2\sin\pi\xi}\,\,\exp\Bigl(\frac{\pi i}{12\tau}+\frac{\xi(\xi-1)\pi i}{2\tau}+\frac{\li(e^{2\pi i\xi})}{2\pi i\tau}+G(\tau,\xi)+P(\tau,\xi)\Bigr)\,. 
\end{align*}

{\bf Remarks --} {\sl {\bf (1)} By using \eqref{equation:Euler} and taking into account \eqref{equation:xqmain}, one can get the Ramanujan's asymptotic expansion formula  \cite[Entry 6, p. 265]{Be1} or its variant \cite[Entry 6', p. 268]{Be1}. 

\medskip
{\bf (2)}
From Theorem \ref{theorem:xqmain} and relations \eqref{equation:symetricG*P}-\eqref{equation:Landen}, one can easily deduce the modular formula for Dedelind $\eta$-function  and that for Jacobi $\theta$-function \cite[\S10.4, p. 496-501]{AAR}. Moreover, when taking the logarithmic derivative with respective to the variable $x$ in \eqref{equation:xqmain}, one can get modular type relations for (generalized) Lambert series.
}

\end{document}